\documentclass{article}
\usepackage{metric-min}

\def\thetitle{Monotonicity of saddle maps}
\def\theauthors{Anton Petrunin and Stephan Stadler}

\hypersetup{colorlinks=true,
citecolor=black,
linkcolor=black,
anchorcolor=black,
filecolor=black,
menucolor=black,
urlcolor=black,
pdftitle={\thetitle},
pdfauthor={\theauthors}
}

\begin{document}

\title{\thetitle}
\author{\theauthors}

\nofootnote{A.~Petrunin was partially supported by NSF grant DMS 1309340.
S.~Stadler was supported by DFG grants STA 1511/1-1 and SPP 2026.}

\newcommand{\Addresses}{{\bigskip\footnotesize
Anton Petrunin, \par\nopagebreak\textsc{Department of Mathematics, PSU, University Park, PA 16802, USA}
\par\nopagebreak
\textit{Email}: \texttt{petrunin@math.psu.edu}

\medskip

Stephan Stadler,
\par\nopagebreak\textsc{Mathematisches Institut der Universit\"at M\"unchen, Theresienstr. 39, D-80333 M\"unchen, Germany}
\par\nopagebreak
\textit{Email}: \texttt{stadler@math.lmu.de}
}}

\date{}

\maketitle

\begin{abstract}
We prove an analog of the Schoen--Yau univalentness theorem for saddle maps between discs.
\end{abstract}

\section{Introduction}

A map from a closed unit disc $\DD$ to a Euclidean space is called \emph{saddle} if for any hyperplane $\Pi$ each connected component of the 
complement $\DD\setminus f^{-1}\Pi$ intersects the boundary~$\partial \DD$.
A connected component that does not meet the boundary will be called a \emph{hat},
so saddle maps can be defined as \emph{maps with no hats}.

Ruled surfaces and harmonic maps provide examples of saddle maps.

We prove a synthetic analog of the univalentness theorem for saddle maps from disc to disc.
The original theorem is formulated for harmonic maps and surfaces with nonpositive curvature.
It was proved by Richard Schoen and Shing Tung Yau \cite{schoen-yau};
an interesting generalization was obtained by J\"urgen Jost \cite{jost}.
An extensive study of general saddle maps was given by Samuil Shefel \cite{shefel-2D,shefel-3D};
his work was inspired by a note of Alexandr Alexandrov \cite{A} on the intrinsic metric of a general ruled surface.
Part of Shefel's work is written in a more popular way \cite{akp}.

A continuous map $f\:X\to Y$ is called \emph{light} if the inverse image of any point $y\in Y$ is totally disconnected.

\begin{thm}{Baby theorem}\label{baby}
Let $f\:\DD\to \DD$ be a light saddle map.
Assume that 
the restriction $f|_{\partial\DD}$ is the identity map.
Then $f$ is a homeomorphism.
\end{thm}

Let $\Delta$ be a surface with possibly nonempty boundary and geodesic metric.
A map $f$ from a closed disc $\DD$ to $\Delta$ is called \emph{saddle} 
if for any geodesic $[x,y]$ in $\Delta$, each connected component of the complement $\DD\setminus f^{-1}[x,y]$ meets the boundary~$\partial\DD$.
It is easy to see that if $\Delta$ is the plane then this new definition of saddle map agrees with the one given above.

A continuous map $f\:X\to Y$ is called \emph{monotone} if the inverse image of any point $y\in Y$ is connected.
Since a connected space is nonempty by definition, \textit{any monotone map is onto}.

\begin{thm}{Monotonicity theorem}\label{thm:main}
Let $\Delta=(\DD,|{*}-{*}|)$ be a closed disc equipped with a $\CAT(\kappa)$ metric 
such that any two points are joined by a unique geodesic, and the geodesic depends continuously on the endpoints.
Assume $f\:\DD\to \Delta$ is a saddle map and the restriction of $f$ to the boundary $\partial\Delta$ is a monotone map 
$\partial\DD\to\partial\Delta$.
Then $f$ is monotone. 
\end{thm}

The theorem essentially states that for saddle maps, monotonicity is an appropriate generalization of univalentness.

Note that the class of spaces satisfying the assumption for $\Delta$ contains any disc with a $\CAT(0)$ metric.
Moreover, this class includes nonpositively curved surfaces considered in the original univalentness theorem \cite{schoen-yau}
as well as the surfaces described in \cite{jost}.

For the saddle map $f\:\DD\to \Delta$ in the theorem,
one can use the so-called \emph{monotone-light factorization} \cite{eilenberg}
which is $f=g\circ h$,
where $g$ is monotone and $h$ is light.
Moore's quotient theorem \cite{moore} implies that the target of the map $g$ is homeomorphic to the disc $\DD$, so we may think that $g\:\DD\to \DD$ and $h\:\DD\to \Delta$.
The light map $h$ is saddle, and according to the theorem, it is monotone.
The latter implies that $h$ is a homeomorphism.

The following theorem is a generalization of the problem ``Saddle surface'' in~\cite{petrunin-orthodox}.
Together with Shefel's theorem (see \cite{shefel-3D} and \cite[4.5.5]{akp}) it implies that the induced 
intrinsic metric on $\Sigma$ satisfies the $\CAT(0)$ comparison. 
(It is not known whether any saddle surface in the Euclidean space has locally $\CAT(0)$ induced intrinsic metric.)

\begin{thm}{Saddle graph theorem}\label{cor:projection}
Let $\Sigma$ be a saddle surface in $\RR^3$ homeomorphic to a disc.
Assume that the orthogonal projection to the $(x,y)$-plane
maps the boundary of $\Sigma$
injectively to a convex closed curve.
Then the orthogonal projection to the $(x,y)$-plane is injective on all of~$\Sigma$.

In particular, $\Sigma$ is a graph $z=f(x,y)$ for a function $f$ defined on a convex figure in the $(x,y)$-plane.
\end{thm}

\medskip

\parbf{Acknowledgments.}
We want to thank Alexander Lytchak and the anonymous referee for keen comments.

\section{Energy-minimizing maps are saddle}

Recall that harmonic maps $f\:M\to N$ between Riemannian manifolds can be defined as local energy minimizers 
among maps with fixed values on the boundary.
Here the energy is defined as 
\[E(f)=\int_M|df|^2,\]
where $df\:\T M\to \T N$ is the differential of~$f$.

\begin{thm}{Proposition} 
Assume $\Delta=(\DD,|{*}-{*}|)$ is a disc with Riemannian metric such that any two points $x,y\in\Delta$ are joined by a minimizing geodesic~$[x,y]$.
Then any energy-minimizing harmonic map $f\:\DD\to\Delta$ with fixed values on the boundary is saddle.
\end{thm}

\parit{Proof.}
Assume the contrary; that is, for some geodesic $[x,y]$ in $\Delta$ the complement $\DD\setminus f^{-1}[x,y]$ has a component 
$\Omega$ which does not meet the boundary~$\partial\DD$.

{

\begin{wrapfigure}{r}{33 mm}
\begin{lpic}[t(-0 mm),b(-0 mm),r(0 mm),l(0 mm)]{pics/cutting-hat(1)}
\lbl{15,17,51;$f(\Omega)$}
\lbl[rb]{9,6;$x$}
\lbl[rb]{25,26;$y$}
\lbl{23,10;$\Delta$}
\end{lpic}
\end{wrapfigure}

Let $\gamma(t)$ be the unit speed parametrization of the geodesic $[x,y]$ from $x$ to~$y$.
Let us redefine the map $f$ in $\Omega$ by setting $t(z)=\min\{\,|x-f(z)|,|x-y|\}$ and
\[\hat f(z)=
\left[
\begin{aligned}
&\gamma\circ t(z)&&\text{if}&& z\in\Omega,
\\
&f(z)&&\text{if}&& z\notin\Omega.
\end{aligned}
\right.\]
Note that $f|_{\partial \DD}\equiv \hat f|_{\partial \DD}$ and
\[|\hat f(z)-\hat f(z')|\le |f(z)-f(z')|\]
for any $z,z'\in \DD$.
It follows that $E(f)>E(\hat f)$ --- a contradiction.
\qeds

}

\parbf{Remark.}
The above proof remains valid if the metric on the target disc is not Riemannian, but has an upper curvature bound in the sense of Alexandrov;
the harmonic maps into such spaces can also be defined as local energy minimizers, see \cite{GS,KS}.

\section{Monotonicity}

\begin{thm}{Claim}\label{claim}
Let $f\:\DD\to \Delta$ be as in the monotonicity theorem (\ref{thm:main}).
Then 
\begin{enumerate}[(i)]
\item\label{claim:i} for any closed convex set $K\subset\Delta$ each connected component of $\DD\setminus f^{-1}K$ intersects~$\partial\DD$.
\item\label{claim:ii} for any open convex set $\Phi\subset\Delta$ each connected component of $f^{-1}\Phi$ is simply connected.
\end{enumerate}
\end{thm}

\parit{Proof; (i).} 
Let $\gamma$ be a geodesic in $\Delta$, disjoint from $K$ and with endpoints on~$\partial\Delta$.
Let $\Sigma_K$ be the set of all such geodesics.
For each $\gamma\in\Sigma_K$ denote by
$H_\gamma$ the component of $\Delta\setminus\gamma$ which contains~$K$.
Since $K$ is closed and convex and $\Delta$ is homeomorphic to a disc, we have
\[K=\bigcap_{\gamma\in\Sigma_K} H_\gamma.\]
In other words, if $x\notin K$, then some geodesic separates $K$ from~$x$. 

By definition of saddle maps, each connected component of $\DD\setminus f^{-1}H_\gamma$ meets the boundary~$\partial\DD$.
Therefore the same holds for the union
\[\DD\setminus f^{-1}K=\bigcup_{\gamma\in\Sigma_K}(\DD\setminus f^{-1}H_\gamma).\]

\parit{(ii)}
Set $\Psi=f^{-1}\Phi\subset\DD$.
Choose a simple closed curve $\gamma\:\SS^1\to\Psi$;
denote by $\Gamma$ the disc bounded by~$\gamma$.
By \textit{(i)}, $f(\Gamma)$ lies in the convex hull of $f(\gamma)$ which will be denoted by~$K$.
Since $\Phi$ is convex, $K\subset \Phi$.
It follows that $\Gamma\subset\Psi$ and therefore $\gamma$ is contractible in~$\Psi$.

Since $\gamma$ is arbitrary, $\Psi$ is simply connected.
\qeds

\parit{Proof of the monotonicity theorem.}
Since $f|_{\partial\DD}\:\partial\DD \to\partial \Delta$ is monotone, it has degree~$\pm1$.
We can assume that the orientations on $\DD$ and $\Delta$ are chosen so that $\deg f|_{\partial\DD}=1$
and therefore $\deg f=1$;
in particular $f$ is onto.

Assume $f$ is not monotone;
that is, there is a point $x\in \Delta$ such that the inverse image $f^{-1}\{x\}$ is not connected.

Given $s\in\partial \Delta$, consider the open set
\[\Phi_s=\set{y\in\Delta}{\measuredangle [x\,^y_s]<\tfrac\pi2}.\]
Note that $\Phi_s$ and its complement are convex in~$\Delta$.
(Here we use that $\Delta$ is homeomorphic to a disc; the analogous statement in higher dimensions does not hold.)
In particular, the relative boundary $\partial_\Delta\Phi_s$ is a geodesic.

Consider the two open subsets $\Theta\subset \DD$ and $\Omega\subset\DD\times\partial\Delta$ defined as
\begin{align*}
\Theta&=\DD\setminus f^{-1}\{x\},
\\
\Omega&=\set{(z,s)\in \DD\times\partial\Delta}{f(z)\in \Phi_s}.
\end{align*}

The projection $\DD\times\partial\Delta\to \DD$ sends $\Omega$ to~$\Theta$.
Note that the induced homomorphism \textit{$\pi_1\Omega\to \pi_1\Theta$ is onto}.
Indeed, for any point $z\in \Theta$ there is $s\in \partial\Delta$ such that $f(z)\in \Phi_s$ or equivalently $\measuredangle [x\,^{f(z)}_s]<\tfrac\pi2$.
Moreover, since $f|_{\partial\DD}$ is monotone, the set of points $s$ satisfying the above condition is an open arc in $\partial\Delta$.
By \cite[Theorem 2]{ungar}, one can fix a continuous map $z\mapsto s_z$ such that $f(z)\in \Phi_{s_z}$ for any $z\in \Theta$.
Then for any loop $\alpha$ in $\Theta$, the loop $\tilde\alpha(t)=(\alpha(t),s_{\alpha(t)})$ is an $f$-lift of $\alpha$ in $\Omega\subset\DD\times\partial\Delta$.

Note that $\Theta\cap \partial\DD$ is connected.
Indeed,
\begin{itemize}
\item If $x\notin \partial\DD$, then $\Theta\cap \partial\DD=\partial \DD$. 
\item If $x\in \partial\DD$, then since $f|_{\partial\DD}$ is monotone, $\Theta\cap \partial\DD$ is an open arc.
\end{itemize}
Since $\{x\}$ is convex, by Claim \ref{claim}\textit{(\ref{claim:i})}, every connected component of $\Theta$ has to intersect~$\partial\DD$.
It follows that $\Theta$ is connected as well.

Consider the restriction of the projection $\DD\times\partial\Delta\to \partial\Delta$ to $\Omega$;
it has fiber $\Psi_s=f^{-1}\Phi_s$ at $s\in\partial\Delta$.
By \ref{claim}, the set $\Psi_s$ is either empty or simply connected.

Indeed, fix $s\in\partial\Delta$ and assume $\Psi_s\ne \emptyset$.
Since $\Phi_s$ is a complement of a closed convex set,
each connected component of $\Psi_s$ must meet $\sigma_s=\Psi_s\cap\SS^1$.
Since $f$ is monotone, $\sigma_s$ is an open arc in~$\SS^1$.
In particular, $\sigma_s$ is connected and so is~$\Psi_s$.
Since $\Phi_s$ is a convex open set, by Claim \ref{claim}\textit{(\ref{claim:ii})}, $\Psi_s$ is simply connected.

Note that $\Phi_s=\emptyset\iff s=x$.
Applying \cite[Theorem 2]{ungar} again, we get the following.
\begin{itemize}
\item If $x\notin\partial\DD$, then from above $\Psi_s$ is simply connected for any~$s$.
Therefore the projection $\Omega\to \partial\Delta$ induces an isomorphism of fundamental groups; that is $\pi_1\Omega=\ZZ$.
\item If $x\in\partial\DD$, then by a similar reason, we have $\pi_1\Omega=0$.
\end{itemize}

Since $f^{-1}\{x\}$ is not connected, it can be divided into two subsets by a curve in~$\Theta$.
\begin{itemize}
\item If $x\notin\partial \DD$, it follows that $\pi_1\Theta$ and therefore $\pi_1\Omega$ contains a free group with two generators --- a contradiction.
\item If $x\in\partial\DD$, it follows that $\pi_1\Theta$ and therefore $\pi_1\Omega$ contains $\ZZ$ as a subgroup --- a contradiction again.\qeds
\end{itemize}

\section{Saddle graph}

The following proof reminds proof of Sergei Bernstein in \cite{bernstein}.

\parit{Proof of the saddle graph theorem (\ref{cor:projection}).}
Denote by $F$ the convex figure in the $(x,y)$-plane bounded by the projection of boundary of~$\Sigma$.
Since $\Sigma$ is saddle, it lies in the convex hull of its boundary;
in particular, the projection of $\Sigma$ lies in~$F$.
Therefore the monotonicity theorem can be applied to the composition of the embedding $\DD\hookrightarrow \RR^3$ which describes $\Sigma$ and the projection
to the $(x,y)$-plane. The composition sends $\DD$ to $F$.

It follows that any vertical line in $F\times \RR$ intersects $\Sigma$ along a closed interval or a point.
In other words, there are two functions $\alpha,\beta\:F\to \RR$
such that 
\[(x,y,z)\in \Sigma
\quad\iff\quad
\alpha(x,y)\le z\le \beta(x,y).\]
Therefore, to prove the saddle graph theorem it is sufficient to show the following:

\begin{thm}{Claim}
$\alpha\equiv \beta$.
\end{thm}

Before proving the claim let us list properties of the functions $\alpha$ and~$\beta$. 

Note that the function $\alpha\:F\to \RR$ is lower semicontinuous
and $\beta\:F\to \RR$ is upper semicontinuous.
Further, $\alpha$ and $\beta$ agree on $\partial F$, and the 
restriction $\alpha|_{\partial F}=\beta|_{\partial F}$ is continuous.
In particular, if an arc $\theta\subset \partial F$ is sufficiently short, then
the change of values of $\alpha$ in $\theta$ is small.

Since $\Sigma$ is saddle, for any linear function $\lambda\:\RR^2\to\RR$ and any closed domain $D\subset F$, 
the restriction $(\alpha-\lambda)|_D$ admits a maximum point on the boundary $\partial D$ --- if not, then a plane $z=\lambda(x,y)+c$ cuts a hat in $D$.

Now let us start to prove the claim.

Assume the contrary; that is, $\alpha(s)<\beta(s)$ for some $s\in F$.
Note that $s$ lies in the interior of $F$ and we can find a pair of points $p$ and $q$ arbitrarily close to $s$
such that 
\[\alpha(q)>\beta(p)+\eps\]
for some fixed $\eps>0$.

Draw the chord $[a,b]$ of $F$ containing $p$ and $q$;
we assume that the points $a,p,q,b$ appear on the chord in this order.

Consider the linear function $\lambda\:[a,b]\to \RR$ such that 
\[\lambda(p)=\beta(p)+\tfrac\eps3,\quad\lambda(q)=\alpha(p)-\tfrac\eps3.\]
Let $\ell$ be the graph of $\lambda$ over $[a,b]$; it is a line segment in $F\times \RR$.

Consider a one-parameter family of planes $\Pi_t\supset \ell$, $t\in(0,1)$ which rotates by angle $\pi$ approaching the vertical plane as $t$ goes to $0$ or $1$.

Note that the plane $\Pi_t$ is a graph of a linear function $\lambda_t\:\RR^2\to \RR$ 
and the restriction $\lambda_t|_{[a,b]}$ agrees with $\lambda$ for any $t$.
Set 
\[\Theta_t=\set{u\in F}{\alpha(u)>\lambda_t(u)};\]
note that $\Theta_t$ is open in $F$ and $\Theta_t\ni a,q$ and $\Theta_t\not\ni p,b$ for any $t$.

We say that the value $t$ is left (right) if $a$ can be connected to $q$ by a curve in $\Theta_t$ which goes on the left (correspondingly right) from $p$.

Since $\Sigma$ is saddle, a value $t$ can not be left and right at the same time --- if it is the case, then $\Pi_t$ cuts a hat from $\Sigma$.

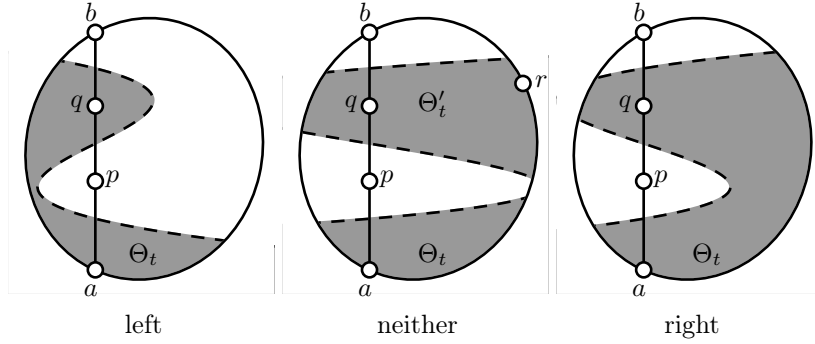
\begin{center}
\begin{lpic}[t(-0 mm),b(5 mm),r(0 mm),l(0 mm)]{pics/yin-yang--left(1)}
\lbl[t]{11,1;$a$}
\lbl[b]{11,36;$b$}
\lbl[l]{13,15;$p$}
\lbl[r]{10,25;$q$}
\lbl{18,5;$\Theta_t$}
\lbl[t]{18,-3;left}
\end{lpic}
\begin{lpic}[t(-0 mm),b(5 mm),r(0 mm),l(0 mm)]{pics/yin-yang--middle(1)}
\lbl[t]{11,1;$a$}
\lbl[b]{11,36;$b$}
\lbl[l]{13,15;$p$}
\lbl[r]{10,25;$q$}
\lbl[l]{33.6,28;$r$}
\lbl{20,25;$\Theta_t'$}
\lbl{20,5;$\Theta_t$}
\lbl[t]{18,-3;neither}
\end{lpic}
\begin{lpic}[t(-0 mm),b(5 mm),r(0 mm),l(0 mm)]{pics/yin-yang--right(1)}
\lbl[t]{11,1;$a$}
\lbl[b]{11,36;$b$}
\lbl[l]{13,15;$p$}
\lbl[r]{10,25;$q$}
\lbl{20,5;$\Theta_t$}
\lbl[t]{18,-3;right}
\end{lpic}
\end{center}

The sets of left (right) values are open in $[0,1]$.
We can assume that all values $t\approx0$ are left and all values $t\approx1$ are right.
Therefore some value, say $t_0$, is neither left nor right.

Note that the connected component $\Theta_{t_0}'\ni q$ in $\Theta_{t_0}$ does not contain $a$.

Since $\Sigma$ is saddle, the restriction $(\alpha-\lambda_{t_0})|_{\Theta_{t_0}'}$ admits its maximum on $\partial F\zz\cap \Theta_{t_0}'$.
By construction $\alpha(q)-\lambda_{t_0}(q)>\tfrac\eps3$; 
therefore 
\[\alpha(r)-\lambda_{t_0}(r)>\tfrac\eps3
\eqlbl{eq:alpha-lambda(r)}
\] 
for some $r\in \partial F\cap \Theta_{t_0}'$.

Denote by $\theta$ the connected component of $r$ in $\partial F\cap \Theta_{t_0}'$.
Without loss of generality, we can assume that the $z$-coordinate of points on $\Sigma$ are between 0 and 1.
It follows that $\theta$ lies in a thin strip $S$ described by the inequalities
\[0\le \lambda_{t_0}\le 1.\]

Since $p$ and $q$ are close to $s$, 
 the slope of $\lambda$ and therefore the slope of $\lambda_{t_0}$ can be made arbitrarily large;
therefore the strip $S$ can be assumed to be arbitrarily thin.
It follows that the arc $\theta$ can be assumed to be  arbitrarily short.
(Here we use that $s$ lies in the interior of $F$ and $S$ passes arbitrarily close to $s$.)
Therefore the change of values of $\alpha$ in $\theta$ can be assumed to be arbitrarily small.

Note that $\lambda_{t_0}$ is monotonic on $\theta$ and 
the function $\alpha-\lambda_{t_0}$ vanishes on the ends of $\theta$.
It follows that 
\[\sup\set{\alpha(u)-\lambda_{t_0}(u)}{u\in\theta}\]
can be assumed to be arbitrarily small.
The latter contradicts \ref{eq:alpha-lambda(r)}.
\qeds

\section{Final remarks}

\parbf{On the baby theorem.}
The presented proof is a tricky fix of the following \emph{fake} proof of the baby theorem (\ref{baby}).
We say where we cheat in the footnote.

\parit{Fake proof.}
Note that $\deg f=1$;
in particular $f$ is onto.
It remains to show that $f$ is injective.

Assume $w=f(x)=f(y)$ for distinct points $x,y\in\DD$.
Note that $w$ lies in the interior of~$\DD$.
Choose a chord $\gamma$ that contains $w$ and goes 
from boundary to boundary of~$\DD$.
The inverse image $p^{-1}(\gamma)$ is a contractible set with two ends at $\partial\DD$, say $a$ and~$b$.
We can assume that the points $a,x,y,b$ appear in the same order on $p^{-1}(\gamma)$.%
\footnote{This is where we are cheating: the inverse image $p^{-1}(\gamma)$ might be as terrible as a pseudoarc, where the order of points has no sense.}

There is a continuous one-parameter family of chords $\gamma_t\ni w$ with the ends at $\partial \DD$
such that $\gamma=\gamma_0$ and $\gamma_1$ is $\gamma$ with reversed parametrization.
Note that the order of $x$ and $y$ on $p^{-1}(\gamma_t)$ does not change in~$t$.
On the other hand, the orders on $\gamma_0$ and $\gamma_1$ are opposite, a contradiction.\qeds

A correct proof of the baby theorem can be built on the deep theorem of Shefel \cite{shefel-2D}, but this argument does not seem to be generalizable.

\parit{Proof.}
Let us extend the map of the disc by the identity map outside the disc. 
According to Shefel's theorem, the induced length metric on the plane is $\CAT(0)$.
This metric coincides with the Euclidean metric outside a compact set 
therefore the induced intrinsic metric is flat and the map preserves this metric;
in particular, it is a homeomorphism.\qeds

{

\begin{wrapfigure}{r}{34 mm}
\begin{lpic}[t(-6 mm),b(-0 mm),r(0 mm),l(0 mm)]{pics/mapping-cylinder(1)}
\end{lpic}
\end{wrapfigure}

\parbf{A generalization of the monotonicity theorem.}
In the proof the condition that $\Delta$ is a disc can be relaxed to the following:
the mapping cylinder of the target space over $f|_{\partial\DD}$ is homeomorphic to a closed disc.
For example, the target space $\Delta$ might look like the solid figure eight in the picture.

Under the name \emph{disc retracts}, these spaces are used in our paper \cite{petrunin-stadler} on a closely related subject. 

}

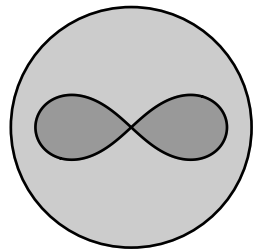
\parbf{On univalentness of harmonic maps.}
If $f\:\DD\to \Sigma$ is a harmonic map from a closed disc to a surface with a Riemannian metric,
then one can show that for any $y\notin f(\partial\DD)$ the inverse image $f^{-1}\{y\}$ is a discrete set of points.
A proof of this statement was suggested by Alexandre Eremenko \cite{eremenko}.
In particular, $f$ is light if $f^{-1}\{f(x)\}=\{x\}$ for any $x\in\partial \DD$.
It reduces the univalentness theorems \cite{schoen-yau,jost} to our monotonicity theorem.

As shown by Ernst Kuwert \cite[Theorem 3]{Ku} the corresponding statement for surfaces with $\CAT(0)$ metrics is wrong. 
Namely, there are (1) a disc $\Delta$ with a flat metric everywhere except for one cone point $q$ where it has negative curvature and (2) a harmonic map $f\:\DD\to \Delta$ such that $f$ restricts to a homeomorphism between boundaries but $f^{-1}\{q\}$ is a nontrivial tree.

Therefore, the monotonicity theorem is optimal even for harmonic maps.

\parbf{On saddle graph theorem.}
The property \ref{claim}\textit{(\ref{claim:i})} can be used to define saddle discs in arbitrary $\CAT(0)$ space.
With this definition, the saddle graph theorem does not hold in the product space $C\times \RR$, where $C$ is a cone with large total angle --- if $o$ is the tip of the cone then the vertical line $\{o\}\times\RR$ might intersect the saddle disc along a nontrivial segment.
An example can be built from the quoted construction of Ernst Kuwert, but actually simpler.

\Addresses

\end{document}